# A Theorem about Simultaneous Orthological and Homological Triangles

Ion Pătraşcu
Fraţii Buzeşti College, Craiova, Romania

Florentin Smarandache
University of New Mexico, Gallup Campus, USA

**Abstract.** In this paper we prove that if $P_1, P_2$ are isogonal points in the triangle $ABC$, and if $A_1B_1C_1$ and $A_2B_2C_2$ are their corresponding pedal triangles such that the triangles $ABC$ and $A_1B_1C_1$ are homological (the lines $AA_1$, $BB_1$, $CC_1$ are concurrent), then the triangles $ABC$ and $A_2B_2C_2$ are also homological.

**Introduction.**

In order for the paper to be self-contained, we recall below the main definitions and theorems needed in solving this theorem.

Also, we introduce the notion of *Orthohomological Triangles*, which means two triangles that are simultaneously orthological and homological.

**Definition 1**

In a triangle $ABC$ the Cevians $AA_1$ and $AA_2$ which are symmetric with respect to the angle's $BAC$ bisector are called isogonal Cevians.

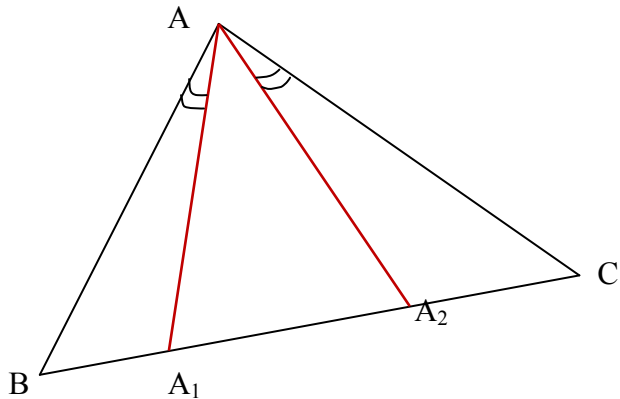

Fig. 1

**Observation 1**

If $A_1, A_2 \in BC$ and $AA_1, AA_2$ are isogonal Cevians then $\sphericalangle BAA_1 \equiv \sphericalangle BAA_2$. (See Fig.1.)

**Theorem 1 (Steiner)**

If in the triangle $ABC$, $AA_1$ and $AA_2$ are isogonal Cevians, $A_1$, $A_2$ are points on $BC$ then:

$$\frac{A_1B}{A_1C} \cdot \frac{A_2B}{A_2C} = \left(\frac{AB}{AC}\right)^2$$

**Proof**
We have:

$$\frac{A_1B}{A_1C} = \frac{\text{area}\Delta BAA_1}{\text{area}\Delta CAA_1} = \frac{\frac{1}{2}AB \cdot AA_1 \sin(\sphericalangle BAA_1)}{\frac{1}{2}AC \cdot AA_1 \sin(\sphericalangle CAA_1)} \quad (1)$$

$$\frac{A_2B}{A_2C} = \frac{\text{area}\Delta BAA_2}{\text{area}\Delta CAA_2} = \frac{\frac{1}{2}AB \cdot AA_2 \sin(\sphericalangle BAA_2)}{\frac{1}{2}AC \cdot AA_2 \sin(\sphericalangle CAA_2)} \quad (2)$$

Because $\sin(\sphericalangle BAA_1) = \sin(\sphericalangle BAA_2)$ and $\sin(\sphericalangle BAA_2) = \sin(\sphericalangle CAA_1)$
by multiplying the relations (1) and (2) side by side we obtain the Steiner relation:

$$\frac{A_1B}{A_1C} \cdot \frac{A_2B}{A_2C} = \left(\frac{AB}{AC}\right)^2 \quad (3)$$

**Theorem 2**
In a given triangle, the isogonal Cevians of the concurrent Cevians are concurrent.
**Proof**
We'll use the Ceva's theorem which states that the triangle's $ABC$ Cevians $AA_1$, $BB_1$, $CC_1$ ( $A_1 \in BC$, $B_1 \in AC$, $C_1 \in AB$ ) are concurrent if and only if the following relation takes place:

$$\frac{A_1B}{A_1C} \cdot \frac{B_1C}{B_1A} \cdot \frac{C_1A}{C_1B} = 1 \quad (4)$$

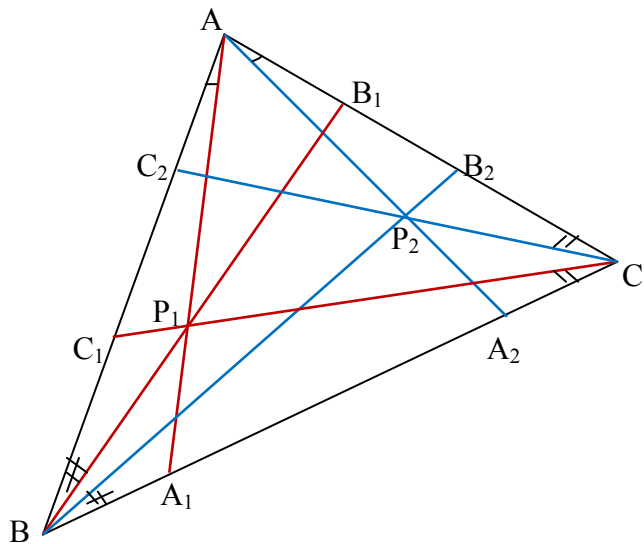

Fig. 2

We suppose that $AA_1$, $BB_1$, $CC_1$ are concurrent Cevians in the point $P_1$ and we'll prove that their isogonal $AA_2$, $BB_2$, $CC_2$ are concurrent in the point $P_2$. (See Fig. 2).

From the relations (3) and (4) we find:

$$\frac{A_2B}{A_2C} = \left(\frac{AB}{AC}\right)^2 \cdot \frac{A_1C}{A_1B} \tag{5}$$

$$\frac{B_2C}{B_2A} = \left(\frac{BC}{AB}\right)^2 \cdot \frac{B_1A}{B_1C} \tag{6}$$

$$\frac{C_2A}{C_2B} = \left(\frac{AC}{BC}\right)^2 \cdot \frac{C_1B}{C_1A} \tag{7}$$

By multiplying side by side the relations (5), (6) and (7) and taking into account the relation (4) we obtain:

$$\frac{A_2B}{A_2C} \cdot \frac{B_2C}{B_2A} \cdot \frac{C_2A}{C_2B} = 1,$$

which along with Ceva's theorem proves the proposed intersection.

**Definition 2**
The intersection point of certain Cevians and the point of intersection of their isogonal Cevians are called isogonal conjugated points or isogonal points.

**Observation 2**
The points $P_1$ and $P_2$ from Fig. 2 are isogonal conjugated points.

In a non right triangle its orthocenter and the circumscribed circle's center are isogonal points.

**Definition 3**
If $P$ is a point in the plane of the triangle $ABC$, which is not on the triangle's circumscribed circle, and $A'$, $B'$, $C'$ are the orthogonal projections of the point $P$ respectively on $BC$, $AC$, and $AB$, we call the triangle $A'B'C'$ the pedal triangle of the point $P$.

**Definition 4**
The pedal triangle of the center of the inscribed circle in the triangle is called the contact triangle of the given triangle.

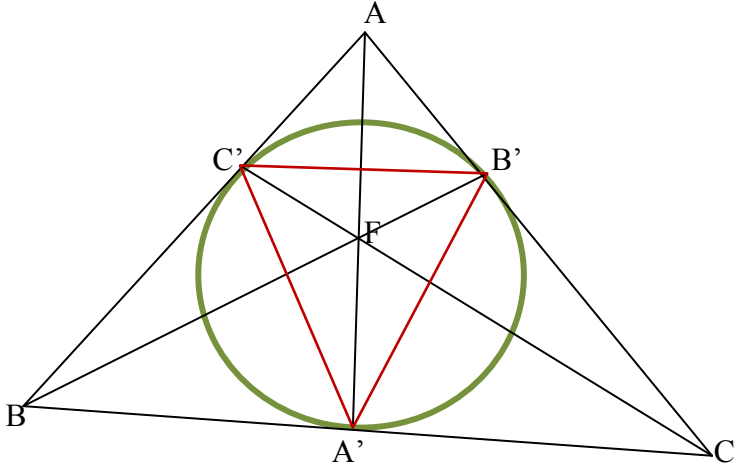

Fig. 3

### Observation 3
In figure 3, $A'B'C'$ is the contact triangle of the triangle $ABC$. The name is connected to the fact that its vertexes are the contact points (of tangency) with the sides of the inscribed circle in the triangle $ABC$.

### Definition 5
The pedal triangle of the orthocenter of a triangle is called orthic triangle.

### Definition 6
Two triangles are called orthological if the perpendiculars constructed from the vertexes of one of the triangle on the sides of the other triangle are concurrent.

### Definition 7
The intersection point of the perpendiculars constructed from the vertexes of a triangle on the sides of another triangle (the triangles being orthological) is called the triangles' orthology center.

### Theorem 3 (The Orthological Triangles Theorem)
If the triangles $ABC$ and $A'B'C'$ are such that the perpendiculars constructed from $A$ on $B'C'$, from $B$ on $A'C'$ and from $C$ on $A'B'$ are concurrent (the triangles $ABC$ and $A'B'C'$ being orthological), then the perpendiculars constructed from $A'$ on $BC$, from $B'$ on $AC$, and from $C'$ on $AB$ are also concurrent.
To prove this theorem firstly will prove the following:

### Lemma 1 (Carnot)
If $ABC$ is a triangle and $A_1$, $B_1$, $C_1$ are points on $BC$, $AC$, $AB$ respectively, then the perpendiculars constructed from $A_1$ on $BC$, from $B_1$ on $AC$ and from $C_1$ on $AB$ are concurrent if and only if the following relation takes place:
$$A_1B^2 - A_1C^2 + B_1C^2 - B_1A^2 + C_1A^2 - C_1B^2 = 0 \qquad (8)$$

**Proof**

If the perpendiculars in $A_1$, $B_1$, $C_1$ are concurrent in the point $M$ (see Fig. 4), then from Pythagoras theorem applied in the formed right triangles we find:

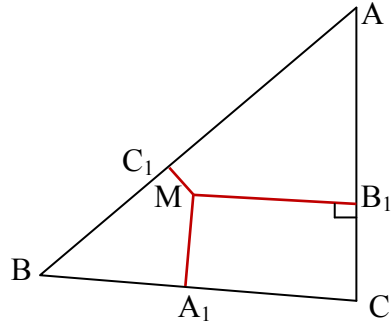

Fig. 4

$$A_1B^2 = MB^2 - MA_1^2 \qquad (9)$$
$$A_1C^2 = MC^2 - MA_1^2 \qquad (10)$$

hence

$$A_1B^2 - A_1C^2 = MB^2 - MC^2 \qquad (11)$$

Similarly it results

$$B_1C^2 - B_1A^2 = MC^2 - MA^2 \qquad (12)$$
$$C_1A^2 - C_1B^2 = MA^2 - MC^2 \qquad (13)$$

By adding these relations side by side it results the relation (8).

**Reciprocally**

We suppose that relation (8) is verified, and let's consider the point $M$ being the intersection of the perpendiculars constructed in $A_1$ on $BC$ and in $B_1$ on $AC$. We also note with $C'$ the projection of $M$ on $AC$. We have that:

$$A_1B^2 - A_1C^2 + B_1C^2 - B_1A^2 + C_1A^2 + C'A^2 - C'B^2 = 0 \qquad (14)$$

Comparing (8) and (14) we find that

$$C_1A^2 - C_1B^2 = C'A^2 - C'B^2$$

and

$$(C_1A - C_1B)(C_1A + C_1B) = (C'A - C'B)(C'A + C'B)$$

and because

$$C_1A - C_1B = C'A + C'B = AB$$

we obtain that $C' = C_1$, therefore the perpendicular in $C_1$ passes through $M$ also.

**Observation 4**

The triangle $ABC$ and the pedal triangle of a point from its plane are orthological triangles.

**The proof of Theorem 3**

Let's consider $ABC$ and $A'B'C'$ two orthological triangles (see Fig. 5). We note with $M$ the intersection of the perpendiculars constructed from $A$ on $B'C'$, from $B$ on $A'C'$ and from $C$ on $A'B'$, also we'll note with $A_1$, $B_1$, $C_1$ the intersections of these perpendiculars with $B'C'$, $A'C'$ and $A'B'$ respectively.

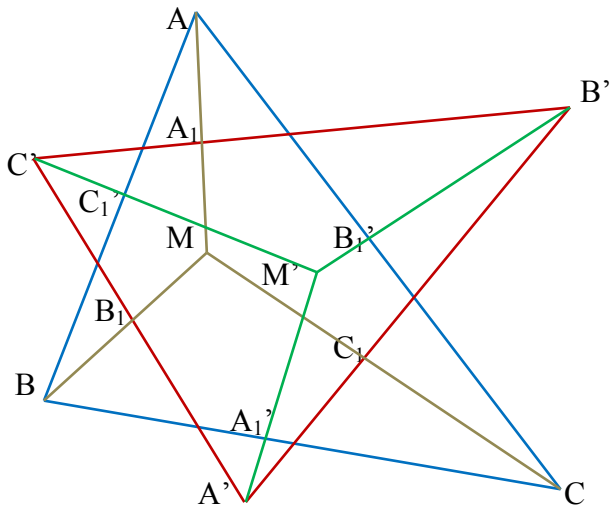

Fig. 5

In conformity with lemma 1, we have:
$$A_1 B'^2 - A_1 C'^2 + B_1 C'^2 - B_1 A'^2 + C_1 A'^2 - C_1 B'^2 = 0 \tag{15}$$

From this relation using the Pythagoras theorem we obtain:
$$B'A^2 - C'A^2 + C'B^2 - A'B^2 + A'C^2 - B'C^2 = 0 \tag{16}$$

We note with $A_1'$, $B_1'$, $C_1'$ the orthogonal projections of $A'$, $B'$, $C'$ respectively on $BC$, $CA$, $AB$. From the Pythagoras theorem and the relation (16) we obtain:

$$A_1'B^2 - A_1'C^2 + B_1'C^2 - B_1'A^2 + C_1'A^2 - C_1'B^2 = 0 \tag{17}$$

This relation along with Lemma 1 shows that the perpendiculars drawn from $A'$ on $BC$, from $B'$ on $AC$ and from $C'$ on $AB$ are concurrent in the point $M'$.

The point $M'$ is also an orthological center of triangles $A'B'C'$ and $ABC$.

**Definition 8**

The triangles $ABC$ and $A'B'C'$ are called bylogical if they are orthological and they have the same orthological center.

**Definition 9**

Two triangles $ABC$ and $A'B'C'$ are called homological if the lines $AA'$, $BB'$, $CC'$ are concurrent. Their intersection point is called the homology point of triangles $ABC$ and $A'B'C'$

**Observation 6**

In figure 6 the triangles $AA'$, $BB'$, $CC'$ are homological and the homology point being $O$

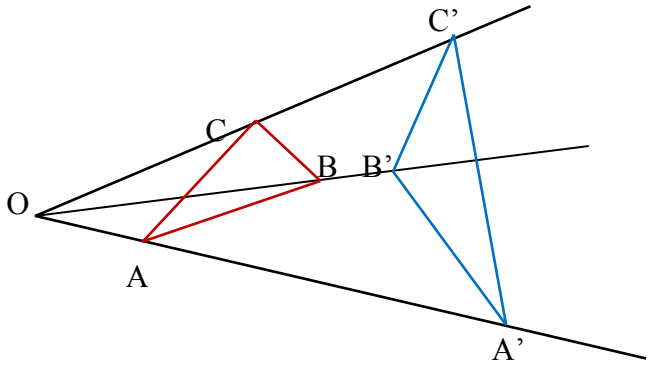

Fig.6

If $ABC$ is a triangle and $A'B'C'$ is its pedal triangle, then the triangles $ABC$ and $A'B'C'$ are homological and the homology center is the orthocenter $H$ of the triangle $ABC$

**Definition 10**

A number of $n$ points ($n \geq 3$) are called concyclic if there exist a circle that contains all of these points.

**Theorem 5 (The circle of 6 points)**

If $ABC$ is a triangle, $P_1, P_2$ are isogonal points on its interior and $A_1 B_1 C_1$ respectively

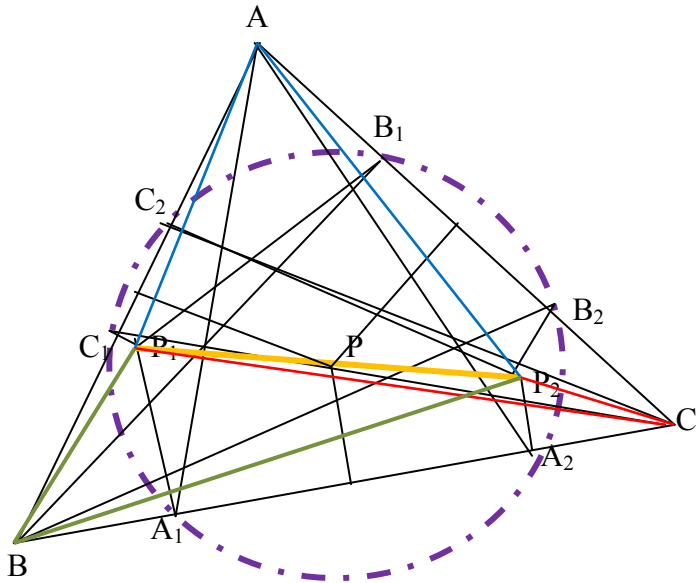

Fig. 7

$A_2B_2C_2$ the pedal triangles of $P_1$ and $P_2$, then the points $A_1, A_2, B_1, B_2, C_1, C_2$ are concyclic.

**Proof**

We will prove that the 6 points are concyclic by showing that these are at the same distance of the middle point P of the line segment $P_1P_2$.

It is obvious that the medians of the segments $(A_1A_2), (B_1B_2), (C_1C_2)$ pass through the point $P$, which is the middle of the segment $(P_1P_2)$. The trapezoid $A_1P_1P_2A_2$ is right angle and the mediator of the segment $(A_1A_2)$ will be the middle line, therefore it will pass through $P$, (see Fig. 7).

Therefore we have:
$$PA_1 = PA_2, \ PB_1 = PB_2, \ PC_1 = PC_2 \tag{18}$$

We'll prove that $PB_1 = PC_2$ by computing the length of these segments using the median's theorem applied in the triangles $P_1B_1P_2$ and $P_1C_2P_2$.

We have:
$$4PB_1^2 = 2(P_1B_1^2 + P_2B_1^2) - P_1P_2^2 \tag{19}$$

We note
$$AP_1 = x_1, \ AP_2 = x_2, \ m(\sphericalangle BAP_1) = m(\sphericalangle CAP_2) = \alpha.$$

In the right triangle $P_2B_2B_1$ applying the Pythagoras theorem we obtain:
$$P_2B_1^2 = P_2B_2^2 + B_1B_2^2 \tag{20}$$

From the right triangle $AB_2P_2$ we obtain:
$$P_2B_2 = AP_2 \sin\alpha = x_2 \sin\alpha \text{ and } AB_2 = x_2 \cos\alpha$$

From the right triangle $AP_1B_1$ it results $AB_1 = AP_1 \cos(A-\alpha)$, therefore
$$AB_1 = x_1 \cos(A-\alpha) \text{ and } P_1B_1 = x_1 \sin(A-\alpha),$$

thus
$$B_1B_2 = |AB_2 - AB_1| = |x_2 \cos\alpha - x_1 \cos(A-\alpha)| \tag{21}$$

Substituting back in relation (17), we obtain:
$$P_2B_1^2 = x_2^2 \sin^2\alpha + [x_2 \cos\alpha - x_1 \cos(A-\alpha)]^2 \tag{22}$$

From the relation (16), it results:
$$4PB_1^2 = 2[x_1^2 + x_2^2 - 2x_1x_2 \cos\alpha \cos(A-\alpha)] P_1P_2^2 \tag{23}$$

The median's theorem in the triangle $P_1C_2P_2$ will give:
$$4PC_2^2 = 2(P_1C_2^2 + P_2C_2^2) - P_1P_2^2 \tag{24}$$

Because $P_1C_1 = x_1 \sin\alpha$, $AC_1 = x_1 \cos\alpha$, $AC_2 = x_2 \cos(A-\alpha)$, $P_1C_2^2 = P_1C_1^2 + C_1C_2^2$, we find that
$$4PC_2^2 = 2[x_1^2 + x_2^2 - 2x_1x_2 \cos\alpha \cos(A-\alpha)] - P_1P_2^2 \tag{25}$$

The relations (23) and (25) show that
$$PB_1 = PC_2 \tag{26}$$

Using the same method we find that :
$$PA_1 = PC_1 \qquad (27)$$
The relations (18), (26) and (27) imply that:
$$PA_1 = PA_2 = PB_1 = PB_2 = PC_1 = PC_2$$
From which we can conclude that $A_1, A_2, B_1, B_2, C_1, C_2$ are concyclic.

**Lemma 2 (The power of an exterior point with respect to a circle)**
If the point $A$ is exterior to circle $C(O,r)$ and $d_1$, $d_2$ are two secants constructed from $A$ that intersect the circle in the points $B, C$ respectively $E, D$, then:
$$AB \cdot AC = AE \cdot AD = cons. \qquad (28)$$
**Proof**
The triangles $ADB$ and $ACE$ are similar triangles (they have each two congruent angles respectively), it results:
$$\frac{AB}{AE} = \frac{AD}{AC}$$

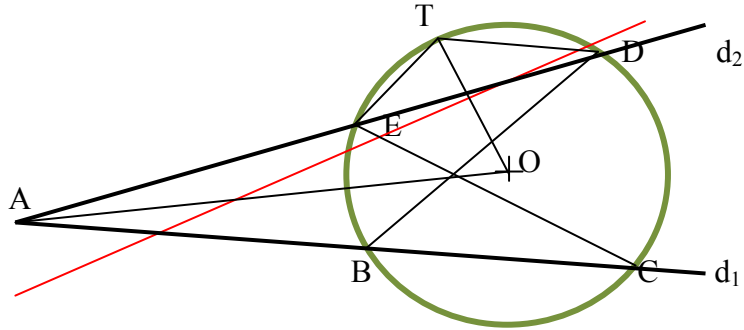

Fig. 8

and from here:
$$AB \cdot AC = AE \cdot AD \qquad (29)$$
We construct the tangent from $A$ to circle $C(O,r)$ (see Fig. 8). The triangles $ATE$ and $ADT$ are similar (the angles from the vertex $A$ are common and $\sphericalangle ATE \equiv \sphericalangle ADT = \frac{1}{2}m(\widehat{TE})$).

We have:
$$\frac{AE}{AT} = \frac{AT}{AD},$$
it results
$$AE \cdot AD = AT^2 \qquad (30)$$
By noting $AO = a$, from the right triangle $ATO$ (the radius is perpendicular on the tangent in the contact point), we find that:
$$AT^2 = AO^2 - OT^2,$$
therefore
$$AT^2 = a^2 - r^2 = const. \qquad (31)$$
The relations (29), (30) and (31) are conducive to relation (28).

**Theorem 6 (Terquem)**

If $AA_1$, $BB_1$, $CC_1$ are concurrent Cevians in the triangle $ABC$ and $A_2$, $B_2$, $C_2$ are intersections of the circle circumscribed to the triangle $A_1, B_1, C_1$ cu $(BC)$, $(CA)$, $(AB)$, then the lines $AA_2$, $BB_2$, $CC_2$ are concurrent.

**Proof**

Let's consider $F_1$ the concurrence point of the Cevians $AA_1$, $BB_1$, $CC_1$.

From Ceva's theorem it results that:
$$A_1B \cdot B_1C \cdot C_1A = A_1C \cdot B_1A \cdot C_1B \tag{32}$$

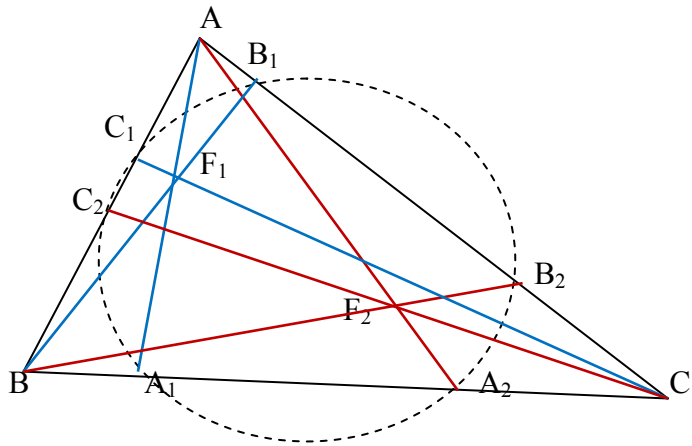

Fig 9

Considering the vertexes $A$, $B$, $C$'s power with respect to the circle circumscribed to the triangle $A_1B_1C_1$, we obtain the following relations:
$$AC_1 \cdot AC_2 = AB_1 \cdot AB_2 \tag{33}$$
$$BA_1 \cdot BA_2 = BC_1 \cdot BC_2 \tag{34}$$
$$CB_1 \cdot CB_2 = CA_1 \cdot CA_2 \tag{35}$$

Multiplying these relations side by side and taking into consideration the relation (32), we obtain
$$AC_2 \cdot BA_2 \cdot CB_2 = AB_2 \cdot BC_2 \cdot CA_2 \tag{36}$$
This relation can be written under the following equivalent format
$$\frac{A_2B}{A_2C} \cdot \frac{B_2C}{B_2A} \cdot \frac{C_2A}{C_2B} = 1 \tag{37}$$

From Ceva's theorem and the relation (37) we obtain that the lines $AA_2$, $BB_2$, $CC_2$ are concurrent in a point noted in figure 9 by $F_2$.

**Note 1**

The points $F_1$ and $F_2$ have been named the Terquem's points by Candido of Pisa – 1900.

For example in a non right triangle the orthocenter $H$ and the center of the circumscribed circle $O$ are Terquem's points.

**Definition 11**
Two triangles are called orthohomological if they are simultaneously orthological and homological.

**Theorem 7**[1]
If $P_1, P_2$ are two conjugated isogonal points in the triangle $ABC$, and $A_1B_1C_1$ and $A_2B_2C_2$ are their respectively pedal triangles such that the triangles $ABC$ and $A_1B_1C_1$ are homological, then the triangles $ABC$ and $A_2B_2C_2$ are also homological.

**Proof**
Let's consider that $F_1$ is the concurrence point of the Cevians $AA_1$, $BB_1$, $CC_1$ (the center of homology of the triangles $ABC$ and $A_1B_1C_1$). In conformity with Theorem 6 the circumscribed circle to triangle $A_1B_1C_1$ intersects the sides $(BC)$, $(CA)$, $(AB)$ in the points $A_2$, $B_2$, $C_2$, these points are exactly the vertexes of the pedal triangle of $P_2$, because if two circles have in common three points, then the two circles coincide; practically, the circle circumscribed to the triangle $A_1B_1C_1$ is the circle of the 6 points (Theorem 5).

Terquem's theorem implies the fact that the triangles $ABC$ and $A_2B_2C_2$ are homological. Their homological center is $F_2$, the second Terquem's point of the triangle $ABC$.

**Observation 7**
If the points $P_1$ and $P_2$ isogonal conjugated in the triangle $ABC$ coincide, then the triangles $ABC$ and $A_2B_2C_2$, the pedal of $P_1 = P_2$ are homological.

**Proof**
From $P_1 = P_2$ and the fact that $P_1, P_2$ are isogonal conjugate, it results that $P_1 = P_2 = I$ - the center of the inscribed circle in the triangle $ABC$. The pedal triangle of $I$ is the contact triangle. In this case the lines $AA_1$, $BB_1$, $CC_1$ are concurrent in $\Gamma$, Gergonne's point, which is the homological center of these triangles.

**Observation 8**
The reciprocal of Theorem 7 for orthohomological triangles is not true.
To prove this will present a counterexample in which the triangle $ABC$ and the pedal triangles $A_1B_1C_1$, $A_2B_2C_2$ of the points $P_1$ and $P_2$ are homological, but the points $P_1$ and $P_2$ are not isogonal conjugated; for this we need several results.

**Definition 12**
In a triangle two points on one of its side and symmetric with respect to its middle are called isometrics.

---
[1] This theorem was called the *Smarandache-Pătraşcu Theorem of Orthohomological Triangles* (see [3], [4]).

**Definition 13**

The circle tangent to a side of a triangle and to the other two sides' extensions of the triangle is called exterior inscribed circle to the triangle.

**Observation 9**

In figure 10 we constructed the extended circle tangent to the side $(BC)$. We note its center with $I_a$. A triangle $ABC$ has, in general, three exinscribed circles

**Definition 14**

The triangle determined by the contact points with the sides (of a triangle) of the exinscribed circle is called the cotangent triangle of the given triangle.

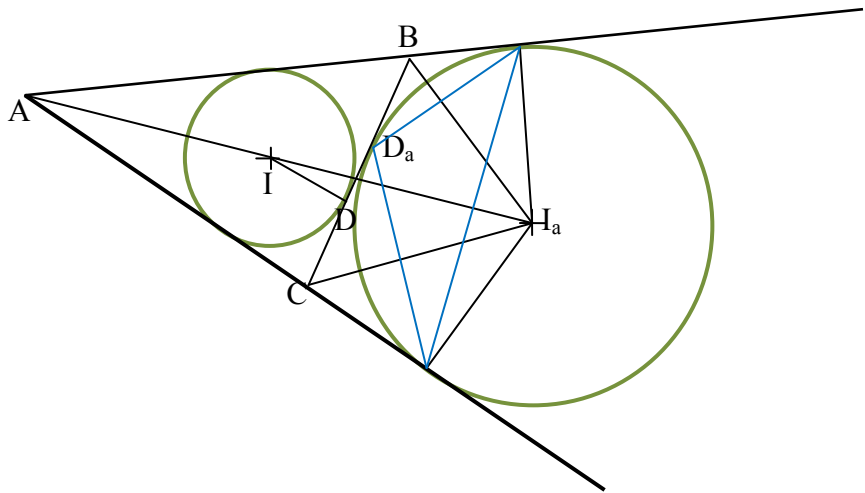

Fig. 10

**Theorem 8**
The isometric Cevians of the concurrent Cevians are concurrent.
The proof of this theorem results from the definition 14 and Ceva's theorem
**Definition 15**
The contact points of the Cevians and of their isometric Cevians are called conjugated isotomic points.

**Lemma 3**
In a triangle $ABC$ the contact points with a side of the inscribed circle and of the exinscribed circle are isotomic points.
**Proof**
The proof of this lemma can be done computational, therefore using the tangents' property constructed from an exterior point to a circle to be equal, we compute the $CD$ and $BD_a$ (see Fig. 10) in function of the length $a,b,c$ of the sides of the triangle $ABC$.

We find that $CD = p - c = BD_a$, which shows that the Cevians $AD$ and $AD_a$ are isogonal ($p$ is the semi-perimeter of triangle $ABC$, $2p = a+b+c$).

**Theorem 9**

The triangle $ABC$ and its cotangent triangle are isogonal.

We'll use theorem 8 and taking into account lemma 3, and the fact that the contact triangle and the triangle $ABC$ are homological, the homological center being the Gergonne's point.

**Observation 10**

The homological center of the triangle $ABC$ and its cotangent triangle is called Nagel's point (N).

**Observation 11**

The Gergonne's point $(\Gamma)$ and Nagel's point (N) are isogonal conjugated points.

**Theorem 10**

The perpendiculars constructed on the sides of a triangle in the vertexes of the cotangent triangle are concurrent.

The proof of this theorem results immediately using lema1 (Carnot)

**Definition 12**

The concurrence point of the perpendiculars constructed in the vertexes of the cotangent triangle on the sides of the given triangle is called the Bevan's point $(V)$.

We will prove now that the reciprocal of the theorem of the orthohomological triangles is false

We consider in a given triangle $ABC$ its contact triangle and also its cotangent triangle. The contact triangle and the triangle $ABC$ are homological, the homology center being the Geronne's point $(\Gamma)$. The given triangle and its cotangent triangle are homological, their homological center being Nagel's point (N). Beven's point and the center of the inscribed circle have as pedal triangles the cotangent triangles and of contact, but these points are not isogonal conjugated (the point $I$ is its own isogonal conjugate).